\documentclass{article}

\usepackage{amsfonts}
\usepackage{graphicx}
\usepackage{amsmath}
\usepackage{fullpage}

\newtheorem{theorem}{Theorem}[section]
\newtheorem{definition}{Definition}[section]

\title{Polynomial Invariant Theory of the Classical Groups}
\author{Quinton Westrich}
\date{March 26, 2006}
\begin{document}
\maketitle

\begin{abstract}
The goal of invariant theory is to find all the generators for the algebra of representations of a group that leave the group invariant. Such generators will be called \emph{basic invariants}. In particular, we set out to find the set of basic invariants for the classical groups GL$(V)$, O$(n)$, and Sp$(n)$ for $n$ even. In the first half of the paper we set up relevant definitions and theorems for our search for the set of basic invariants, starting with linear algebraic groups and then discussing associative algebras. We then state and prove a monumental theorem that will allow us to proceed with hope: it says that the set of basic invariants is finite if $G$ is reductive. Finally we state without proof the First Fundamental Theorems, which aim to list explicitly the relevant sets of basic invariants, for the classical groups above. We end by commenting on some applications of invariant theory, on the history of its development, and stating a useful theorem in the appendix whose proof lies beyond the scope of this work.

\end{abstract}

\section{Linear Algebraic Groups and their Representations}

\subsection{Linear Algebraic Groups}

First we shall establish some notations. We denote by GL$(n,\mathbb{C})$ the group of invertible $n\times n$ complex matrices. By $M_n(\mathbb{C})$ we mean the space of all $n\times n$ complex matrices. We write the $i,j$ entry of a matrix $y \in M_n(\mathbb{C})$ with $1 \leq i,j \leq n$ as $x_{ij}(y)$. The ring of complex polynomials in $n$ variables will be denoted $\mathbb{C}[x_1,\ldots,x_n]$. Now we may state the following definitions:

\begin{definition}
A function $f:M_n(\mathbb{C})\to \mathbb{C}$ is a \textbf{\emph{polynomial function}} on $M_n(\mathbb{C})$ if there exists $p \in \mathbb{C}[x_{11}(y),x_{12}(y),\ldots,x_{nn}(y)]$ such that 
	\[ f(y) = p(x_{11}(y),x_{12}(y),\ldots,x_{nn}(y)) .\]
\end{definition}

\begin{definition}
A subgroup $G \leq \mbox{\emph{GL}}(n,\mathbb{C})$ is a \textbf{\emph{linear algebraic group}} if there exists a set $A$ of polynomial functions on $M_n(\mathbb{C})$ such that 
	\[ G = \{g \in G ~|~ f(g)=0 \mbox{\emph{ for all }} f \in A\} .\]
\end{definition}

A consequence of the \emph{Hilbert basis theorem} (App.\ref{sec:hilbert}) is that any linear algebraic group can be defined by a \emph{finite} number of polynomial equations \cite{goodman}. As an example of a linear algebraic group, consider the subgroup $D_n \leq$GL$(n,\mathbb{C})$ of diagonal invertible matrices. The defining equations for $D_n$ are $x_{ij}=0$ for $i \neq j$. So $D_n$ is a linear algebraic group.

\subsection{Regular Functions}

\begin{definition}
The ring of \textbf{\emph{regular functions}} for the linear algebraic group \emph{GL}$(n,\mathbb{C})$ is defined as
	\[ \mbox{\emph{Aff(GL}}(n,\mathbb{C})):= \mathbb{C}[x_{11}(y),x_{12}(y),\ldots,x_{nn}(y), 							(\det(y))^{-1}]. \]
\end{definition}	
This is just the complex algebra generated by the matrix entry functions $x_{ij}$ and the function (det$(y)$)$^{-1}$.

We must now state a few more notational conventions. We shall denote a vector space by $V$, the algebra of all linear transformations on $V$ by End($V$), and the group of all invertible linear transformations on $V$ by GL($V$).

\begin{definition}
	Suppose $\dim V = n$ and $\{e_i\}$ is a basis for $V$. Let $g \in$ \nolinebreak[3] 					\emph{GL($V$)} and $\phi(g)$ be the matrix $[g_{ij}]$ such that 
		\[ ge_j = \sum_{i=1}^n g_{ij}e_i. \]
	Then the map $g \mapsto \phi (g)$ gives an isomorphism
		\[ \phi : \mbox{\emph{GL}}(V) \stackrel{\cong}{\longrightarrow} \mbox{\emph{GL}} 						(n,\mathbb{C}). \]
	We define the \textbf{\emph{regular functions}} on \emph{GL}$(V)$ to be those of the form 	$f \circ \phi$, where $f$ is a regular function on \emph{GL}$(n,\mathbb{C})$.
\end{definition}
This definition just says that it only makes sense to move into a basis before we can talk about a regular function on GL$(V)$. To denote the algebra of regular functions on GL$(V)$, we use the symbol Aff(GL($V$)). Note, however, that Aff(GL($V$)) is independent on the particular choice of basis of $V$.

\begin{definition}
	Let $G \leq \mbox{\emph{GL}}(V)$ be a linear algebraic group. Let 													$f:\nolinebreak[4]\mbox{\emph{GL}}(V)\to\mathbb{C}$ be a regular function on 								$\mbox{\emph{GL}}(V)$. Then the restriction $f|_G:G\to\mathbb{C}$ is called a 							\textbf{\emph{regular function}} on $G$.
\end{definition}
The set Aff($G$) of regular functions on $G$ is a commutative algebra over $\mathbb{C}$ under pointwise multiplication.

\begin{definition}
	Let $G$ and $H$ both be linear algebraic groups. Then a  \textbf{\emph{regular 							homomorphism}} is a group homomorphism $\phi:G\to H$ such that 															$\phi^*(\mbox{\emph{Aff}}(H))\subseteq \mbox{\emph{Aff}}(G)$, where for $f \in 							\mbox{\emph{Aff}}(H)$, 
		\[\phi^*(f)(g):=f(\phi(g)).\]
\end{definition}
\begin{definition}
	Let $G$ and $H$ both be linear algebraic groups. We say that $G$ and $H$ are 								\textbf{\emph{isomorphic}} as linear algebraic groups if there exists a regular 						homomorphism $\phi:G\to H$ such that $\phi^{-1}$ is also a regular homomorphism.
\end{definition}

\subsection{Representations of Linear Algebraic Groups}

\begin{definition}
	A \textbf{\emph{representation}} of a linear algebraic group $G$ is a pair $(\rho,V)$, 			where $V$ is a complex finite- or infinite-dimensional vector space and 										$\rho:\nolinebreak[4]G\to$\emph{GL}$(V)$ is a group homomorphism.
\end{definition}
\begin{definition}
	A \textbf{\emph{regular representation}} of a linear algebraic group $G$ is a 							representation $(\rho,V)$ for which $\dim (V)$ is finite and $\rho:G\to$\emph{GL}$(V)$ is 	a regular homomorphism.
\end{definition}

This means that if we fix a basis and write out the matrix for $\rho(g)$ in this basis, then the elements $\rho_{ij}(g)$ are all regular functions on $G$.

\begin{definition}
	If $(\rho,V)$ is a regular representation of the linear algebraic group $G$ and $W 					\subseteq V$ is a vector subspace of $V$, then we say that $W$ is 													$G$\textbf{\emph{-invariant}} if $\rho(g)w\in W$ for all $g \in G$, $w \in W$.
\end{definition}
\begin{definition}
	A representation $(\rho,V)$ with $V\neq\{0\}$ is \textbf{\emph{reducible}} if there is a 		non-trivial $G$-invariant subspace $W \subseteq V$. A representation is 										\textbf{\emph{irreducible}} if it is not reducible.
\end{definition}

As an example, let $(\rho,V)$ and $(\sigma,W)$ be regular representations of a linear algebraic group $G$. Define the \emph{tensor product representation} $\rho \otimes \sigma$ on $V \otimes W$ by
	\[(\rho \otimes \sigma)(g)(v \otimes w):= \rho(g)v \otimes \sigma(g)w \]
for $g \in G$, $v \in V$, and $w \in W$. Then $\rho \otimes \sigma$ is a regular representation.

\begin{definition}
	A regular representation $(\rho,V)$ of a linear algebraic group $G$ is 											\textbf{\emph{completely reducible}} if for every $G$-invariant subspace $W \subseteq V$ 		there exists another $G$-invariant subspace $U \subseteq V$ such that $V=W \oplus U$.
\end{definition}

In terms of matrices, this means that any basis $\{w_1,\ldots,w_p\}$ for $W$ can be completed to a basis $\{w_1,\ldots,w_p,u_1,\ldots,u_q\}$ for $V$ so that the subspace $U= \mbox{span}_{\mathbb{C}}\{u_1,\ldots,u_q\}$ is invariant under $\rho(g)$ for all $g\in G$. Thus, the matrix of $\rho(g)$ relative to this basis has the block-diagonal form
\begin{displaymath}	
	\left(
	\begin{array}{cc}
	\sigma (g) & 0 \\
	0				   & \tau(g)
	\end{array}
	\right)
\end{displaymath}
where $\sigma(g)=\rho(g)|_W$ and $\tau(g)=\rho(g)|_U$.

\begin{definition}
	We say that a linear algebraic group $G$ is \textbf{\emph{reductive}} if every regular 			representation $(\rho,V)$ of $G$ is completely reducible.
\end{definition}

In particular, it can be shown that all finite groups and the classical groups are reductive \cite{goodman}.

\section{Associative Algebras and their Representations}

In this section we introduce some of the basic terminology to describe the space $\mathcal{P}(V)$ of polynomials of elements of a vector space $V$. To understand what one means by a polynomial function on a vector space we need the concept of an \emph{algebra}.

\subsection{Associative Algebras}

\begin{definition}
	An \textbf{\emph{associative algebra}} over the complex field $\mathbb{C}$ is a pair 				$(\mathcal{A},\mu)$ with $\mathcal{A}$ a vector space over $\mathbb{C}$ and $\mu$ a 				bilinear associative map on $\mathcal{A}$, i.e.
		\begin{displaymath}
			\mu : \mathcal{A} \times \mathcal{A} \to \mathcal{A}
		\end{displaymath}
	with
		\begin{displaymath}
			(x,y) \mapsto \mu(x,y)\equiv xy 
		\end{displaymath}
	such that
		\begin{displaymath}
			(xy)z=x(yz) 
		\end{displaymath}
	for all $x,y,z \in \mathcal{A}$. If the algebra possesses a unit element, it is called a 		\textbf{\emph{unital algebra}}.
\end{definition}

When speaking of an algebra $(\mathcal{A},\mu)$, it is customary to only write the vector space $\mathcal{A}$ and omit mentioning explicitly the bilinear map $\mu$. Also, we note that it is easy to see that for a vector space $V$ the space of linear maps End$(V)$ is an associative algebra with composition taken as the algebra multiplication.

\begin{definition}
	Let $\mathcal{J}$ be a vector subspace of the algebra $\mathcal{A}$. Then if $\mu 					|_{\mathcal{J}}$ is closed in $\mathcal{J}$, we call $\mathcal{J}$ a 												\emph{\textbf{subalgebra}} of $\mathcal{A}$ with respect to the bilinear map 								$\mu|_{\mathcal{J}}$.
\end{definition}

\begin{definition}
	Let $\mathcal{A}$ be an associative algebra and $\mathcal{I}$ be a subalgebra of 						$\mathcal{A}$. Then $\mathcal{I}$ is called a \emph{\textbf{two-sided ideal}} of the 
	algebra $\mathcal{A}$ if, $\forall (x \in \mathcal{I}), \, \forall (y \in \mathcal{A})$,
	$xy \in \mathcal{I}$ and $yx \in \mathcal{I}$.
\end{definition}

\begin{definition}
	Let $\mathcal{A}$ and $\mathcal{C}$ be associative algebras. An \emph{\textbf{algebra 			homomorphism}} is a map $\varphi:\mathcal{A} \to \mathcal{C}$ such that for all 						$x,y\in\mathcal{A}$,
	\begin{equation}
		\varphi(xy)=\varphi(x)\varphi(y).
	\end{equation}
\end{definition}

\begin{definition}
	An \emph{\textbf{algebra isomorphism}} is a bijective algebra homomorphism.
\end{definition}

\subsection{Representations of Associative Algebras}

\begin{definition}
	A \textbf{\emph{representation}} of an associative algebra $\mathcal{A}$ is a pair 					$(\rho,V)$ with $V$ a vector space and $\rho:\mathcal{A} \to $\emph{End}$(V)$ an algebra 		homomorphism. Here we call $V$ an $\mathcal{A}$\textbf{\emph{-module}}.
\end{definition}

Let $\mathcal{A}$ be an associative algebra and $U$ be a finite-dimensional irreducible $\mathcal{A}$-module. Denote by $[U]$ the equivalence class of all $\mathcal{A}$-modules isomorphic to $U$. Denote by $\widehat{\mathcal{A}}$ the set of all equivalence classes of finite-dimensional irreducible $\mathcal{A}$-modules.

\begin{definition}
	Let $\mathcal{A}$ be an associative algebra and $V$ be a completely reducible 							$\mathcal{A}$-module. Then for $\xi\!\in\!\widehat{\mathcal{A}}$, we define the 						$\xi$\textbf{\emph{-isotypic subspace}} $V_{(\xi)}$ of $V$ to be
		\begin{equation}
			V_{(\xi)}:=\sum_{\substack{ U \subseteq V \\ [U]= \xi }}U
		\end{equation}
	where the $U$ are invariant and irreducible.
\end{definition}

Let $V$ be a completely reducible $\mathcal{A}$-module and
\begin{equation}
	V=\bigoplus_{i=1}^d V_i
\end{equation}
be any decomposition such that all the $V_i$ are invariant and irreducible. Then it can be shown that $\forall (\xi \!\in\! \widehat{\mathcal{A}}\,)$,
\begin{equation}
	V_{(\xi)}=\bigoplus_{[V_j]=\xi}V_j
\end{equation}
and so
\begin{equation}
	V=\bigoplus_{\xi \in \widehat{\mathcal{A}}} V_{(\xi)}.
\end{equation}
We have now developed a sufficient amount of theory to begin our development of invariant theory.

\section{The Ring of Polynomial Invariants}

Let $G$ be a reductive linear algebraic group and $(\pi,V)$ be a regular representation of $G$. Define a representation $\rho$ of $G$ on the algebra $\mathcal{P}(V)$ by
\begin{equation}
	\rho(g)f(v):=f(\pi(g^{-1})v)\equiv f(g^{-1}v)
\end{equation}
for $f \!\in\! \mathcal{P}(V)$, $g \!\in\! G$, and $v \!\in\! V$, where in the last equality we have suppressed the $\pi$ notation as convention dictates. Note that $\mathcal{P}(V)$ contains polynomials of all orders. We'll need to be more specific so we introduce the notation:
\begin{equation}
	\mathcal{P}^k(V):= \{ f \in \mathcal{P}(V) \,|\, f(zv)=z^kf(v) \mbox{ for } z \in 					\mathbb{C}^{\times} \},
\end{equation}
the (finite-dimensional) space of \emph{homogeneous polynomials of degree $k$} for $k \!\in\! \mathbb{N}$. In particular, $\mathcal{P}^k(V)$ is $G$-invariant and $\rho|_{\mathcal{P}^k(V)}\equiv\rho_k$ is a regular representation of $G$. 
\begin{definition}
	The \emph{\textbf{algebra of $G$-invariants}} is the space of $G$-invariant polynomials 		on a vector space $V$. We shall denote this algebra by $\mathcal{P}(V)^G$.
\end{definition}
\begin{theorem} \label{thm:fingen}
	Suppose $G$ is a reductive linear algebraic group which has a regular representation on a 	vector space $V$. Then the algebra $\mathcal{P}(V)^G$ of $G$-invariant polynomials on $V$ 	is finitely generated as an algebra over $\mathbb{C}$.
\end{theorem}
\begin{subparagraph}{Proof.}
Since $G$ is reductive we have that $\mathbb{C}[\rho_k(G)]$ is semisimple. Hence, we can write
\begin{equation} \label{eq:decomp}
	\mathcal{P}^k(V)=\bigoplus_{\sigma \in \widehat{G}}W_{(\sigma)},
\end{equation}
the decomposition into $G$-isotypic subspaces.

Consider some polynomial function $f \in \mathcal{P}(V)$ of arbitrary order. We can decompose $f$ as
\begin{equation}
	f=\sum_{k=0}^d f_k
\end{equation}
with the $f_k$ homogeneous of degree $k$ so that we have $d$ polynomials each in some $\mathcal{P}^k(V)$. Decompose further every $f_k$ by (\ref{eq:decomp}), collect like $\sigma$'s, and write
\begin{equation}
	f=\sum_{\sigma \in \widehat{G}}f_{(\sigma)}
\end{equation}
where $f_{(\sigma)} \in W_{(\sigma)}$ is the $\sigma$\emph{-isotypic component of} $f$. 

Denote by $f_{(1)}=f^{\natural}$ the trivial representation. Let $\varphi$ be a $G$-invariant function, i.e. $\varphi \in \mathcal{P}(V)^G$ and $f \in \mathcal{P}(V)$ as above. Then, since multiplication by a $G$-invariant function leaves isotypic subspaces invariant,
\begin{equation}
	(\varphi f)^{\natural}=\varphi f^{\natural}.
\end{equation}
So, taking $\mathcal{P}(V)$ as a module of $\mathcal{P}(V)^G$, the $\mathcal{P}(V)^G$-module map $f \mapsto f^{\natural}$ is a projection operator.

Now, by the Hilbert basis theorem (App.\ref{sec:hilbert}), every ideal of $\mathcal{P}(V)$ is finitely generated. Also if $\mathcal{I}$ is an ideal of $\mathcal{P}(V)$, then $\mathcal{P}(V)/\mathcal{I}$ is finitely generated too. Denote by $\mathcal{P}(V)_+^G$ the space of invariant polynomials in which the zeroth order term vanishes. We claim that $\mathcal{P}(V)_+^G$ is an ideal of $\mathcal{P}(V)$. Indeed, let $x \in \mathcal{P}(V)_+^G$ and $y \in \mathcal{P}(V)$. Then we can write $y$ as $y=w+z$ with $\mathrm{deg}\,w\geq 1$ and $\mathrm{deg}\,z=0$, i.e. $w \in \mathcal{P}(V)_+^G$ and $z \in \mathbb{C}$. Then 
\begin{equation}
	xy=x (w+z)=x w + x z.
\end{equation}
Obviously, $\mathrm{deg}\,(x w)\geq 1$. But we also have $\mathrm{deg}\,(x z)\geq 1$ since for any polynomials $a,b\in\mathcal{P}(V)$, $\mathrm{deg}\,(ab)=\mathrm{deg}\,(a)+\mathrm{deg}\,(b)$. The proof that $yx \in \mathcal{P}(V)_+^G$ is completely symmetric. It follows that $\mathcal{P}(V)_+^G$ is finitely generated.

Suppose the set $\{\phi_i\}_{i=1}^n$ generates $\mathcal{P}(V)_+^G$. We claim now that $\{\phi_i\}_{i=1}^n$ must also generate $\mathcal{P}(V)^G$ as an algebra over $\mathbb{C}$. Let $\phi \in \mathcal{P}(V)^G$. Then there exists a set of polynomials $\{f_i\}_{i=1}^n$ with $f_i \in \mathcal{P}(V)$ for every $i$ such that
\begin{equation}
	\phi = \sum_{i=1}^n f_i\phi_i.
\end{equation}
If we now project out $\phi^{\natural}$ so that
\begin{equation}
	\phi^{\natural}=\sum_{i=1}^n (f_i\phi_i)^{\natural}=\sum_{i=1}^n f_i^{\natural}\phi_i,
\end{equation}
we may assume, since $\mathrm{deg}\,f_i^{\natural}\leq\mathrm{deg}\,f_i\leq\mathrm{deg}\,\phi$,
that $f_i^{\natural}$ is in the algebra generated by the set $\{\phi_i\}_{i=1}^n$, and, hence, so is $\phi$.$\,_{\square}$
\end{subparagraph}

Thm. \ref{thm:fingen} motivates the following definition.

\begin{definition}
	The smallest set $\{f_i\}_{i=1}^n$ that generates $\mathcal{P}(V)^G$ is called the set of   \emph{\textbf{basic invariants}}.
\end{definition}

The goal of invariant theory is to find this set. In the following we shall state without proof the sets of basic invariants for the classical groups. It should be noted that while the above theorem states the necessary existence of such a set for reductive groups, these sets are not unique in general. However if we order the set of degrees of each of the generators, that set is unique \cite{goodman}.

\section{Polynomial Invariants of the Classical Groups}
\subsection{Polynomial Invariants of GL$(V)$}

Given a vector space $V$, there exist \emph{natural isomorphisms}
\begin{equation}
	(V^*)^{\oplus k} \cong	\mbox{Hom}(V,\mathbb{C}^k)
\end{equation}
given by
\begin{displaymath}
	v \mapsto \left[ \left\langle v_1^*,v \right\rangle, \ldots, \left\langle v_k^*,v 					\right\rangle \right]
\end{displaymath}
and
\begin{equation}
	V^{\oplus m} \cong	\mbox{Hom}(\mathbb{C}^m,V)
\end{equation}
given by
\begin{displaymath}
	\left[ c_1, \ldots, c_m \right] \mapsto c_1v_1+\cdots+c_mv_m.
\end{displaymath}
From this, we can write the \emph{algebra} isomorphism
\begin{equation}
	\mathcal{P}((V^*)^{\oplus k}\times V^{\oplus m}) \cong \mathcal{P} (\mbox{Hom} 							(V,\mathbb{C}^k) \times \mbox{Hom}(\mathbb{C}^m,V))
\end{equation}
Suppose $f \in \mathcal{P}(\mbox{Hom}(V,\mathbb{C}^k) \times \mbox{Hom}(\mathbb{C}^m,V))$ and $g \in$GL$(V)$. We define the \emph{action of $g$ on $f$} by
\begin{equation}
	g \cdot f(x,y) = f(x\rho(g^{-1}),\rho(g)y).
\end{equation}
Denoting by $M_{k,m}$ the (vector) space of all $k \times m$ complex matrices we can also define the function
\begin{displaymath}
	\mu : \mbox{Hom}(V,\mathbb{C}^k) \times \mbox{Hom}(\mathbb{C}^m,V) \to M_{k,m}
\end{displaymath}
by
\begin{displaymath}
	\mu (x,y) := xy,
\end{displaymath}
in which by juxtaposition we mean ordinary composition of linear maps. Then $\forall (g\in$GL$(V)$),
\begin{eqnarray}
	g \cdot \mu(x,y) & = & \mu(x\rho(g^{-1}),\rho(g)y)  \nonumber \\
									& = & x\rho(g^{-1})\rho(g)y         					\\
									& = & xy                            \nonumber \\
									& = & \mu(x,y).											\nonumber
\end{eqnarray}
Now, define
\begin{displaymath}
	\mu^*:\mathcal{P}(M_{k,m})\to\mathcal{P}(\mbox{Hom}(V,\mathbb{C}^k)\times 									\mbox{Hom}(\mathbb{C}^m,V))^{\mathrm{GL}(V)}
\end{displaymath}
such that $\forall (f \in \mathcal{P}(M_{k,m})), \,	f \mapsto f \circ \mu$.
Finally, define the \emph{matrix entry function} on $M_{k,m}$:
\begin{equation}
	z_{ij}:=\mu^*(x_{ij}).
\end{equation}
Then $z_{ij}$ is the \emph{contraction} of the $i$th dual vector and the $j$th vector position, i.e.
\begin{equation}
	z_{ij}(v_1^*,\ldots,v_k^*,v_1,\ldots,v_m)=\left\langle v_i^*,v_j \right\rangle.
\end{equation}
\begin{theorem}[Polynomial FFT for $\mathrm{GL}(V)$]
The map $\mu^*$ is surjective and, hence, $\mathcal{P}((V^*)^{\oplus k}\times
V^{\oplus 		m})^{\mathrm{GL}(V)}$ is generated by the contractions
	\begin{displaymath}
		\{ \left\langle v_i^*,v_j \right\rangle |\, i=1,\ldots,k;j=1,\ldots,m \}.
  \end{displaymath}
\end{theorem}

\subsection{Polynomial Invariants of O$(n)$}

Let $V=\mathbb{C}^n$ and define a nondegenerate symmetric bilinear form on on $V$ by
\begin{equation}
	(x,y):= \sum_{i=1}^n x_iy_i
\end{equation}
for $x,y \in \mathbb{C}^n$. Denote by O$(n)$ the \emph{orthogonal group} for $(\cdot\,,\cdot)$ so that
\begin{equation}
	g \in \mbox{O}(n) \, :\Leftrightarrow \, g^Tg=\mathbf{1}_{n\times n}.
\end{equation}
Denote also the space of all $k \times k$ symmetric matrices over $\mathbb{C}$ by $SM_k$ so that
\begin{displaymath}
	B \in SM_k \Rightarrow B=B^T.
\end{displaymath}
Define the map $\tau : M_{n,k} \to SM_k$ such that $X \mapsto X^TX$ for all $X \in M_{n,k}$.
Then for all $g \in$O$(n)$ and $X \in M_{n,k}$,
\begin{eqnarray}
	\tau(gX) &=& (gX)^TgX \nonumber \\
					&=& X^Tg^TgX 						\\
					&=& X^TX    \nonumber \\
					&=& \tau(X). \nonumber 
\end{eqnarray}
For $f\in\mathcal{P}(SM_k)$, define $\tau^*:\mathcal{P}(SM_k) \to \mathcal{P}(V^k)^{\mathrm{O}(n)}$ by $f \mapsto f \circ \tau$. From this definition, it follows that
\begin{equation} \label{eq:taustinv} 
	\tau^*(f)(gX)=\tau^*(f)(X)
\end{equation}
and $\tau$ is an algebra isomorphism. For example, if $v_1,\ldots,v_k\in\mathbb{C}^n$, then there exists $X=[v_1,\ldots,v_k]\in M_{n,k}$ so that
\begin{equation} \label{eq:sympostocontr} 
	x_{ij}(X^TX)=(v_i,v_j).
\end{equation}
Now, by Eqs.(\ref{eq:taustinv}) and (\ref{eq:sympostocontr}), we have on $(\mathbb{C}^n)^{\oplus m}$:
\begin{equation}
	\tau^*(x_{ij})(v_1,\ldots,v_k)=(v_i,v_j),
\end{equation}
the contraction of the $i,j$th vector using $(\cdot\,,\cdot)$.
\begin{theorem}[Polynomial FFT for $\mathrm{O}(n)$]
The map $\tau^*\in \mathrm{Hom} (\mathcal{P}(SM_k), 				
					\mathcal{P}((\mathbb{C}^n)^{\oplus
k})^{\mathrm{O}(n)})$ is surjective, and, hence, 				
$\mathcal{P}((\mathbb{C}^n)^{\oplus k})^{\mathrm{O}(n)}$ is generated by the 	
						orthogonal contractions
	\begin{displaymath}
		\{(v_i,v_j)\,|\,1\leq i\leq j \leq k \}.
	\end{displaymath}
\end{theorem}

\subsection{Polynomial Invariants of Sp$(n)$}

Define the matrices
\begin{displaymath}
	\kappa := \left(\begin{array}{cc} 0 & 1 \\ -1 & 0 \end{array}\right)
	\quad \mathrm{ and } \quad
	J_n:= \left(\begin{array}{cccc} 
		\kappa & 0           & \cdots      & 0       \\
		0      & \kappa      & \phantom{0} & \vdots  \\
		\vdots & \phantom{0} & \ddots      & 0       \\
		0      & \cdots      & 0           & \kappa
	\end{array} \right)
\end{displaymath}
(where $J_n$ is block-diagonal) and the antisymmetric form
\begin{equation}
	\omega(x,y):=(x,J_ny)
\end{equation}
for $x,y \in \mathbb{C}^n$. Obviously, if we mean by the subscript $n$ the dimension of $J_n$, then $n \in 2\mathbb{N}$. Denote by Sp$(n)$ the \emph{symplectic group} for $\omega(\cdot\,,\cdot)$ so that
\begin{equation}
	g \in \mbox{Sp}(n) \, :\Leftrightarrow \, g^TJ_ng=J_n.
\end{equation}
For the space of $k \times k$ antisymmetric matrices we write $AM_k$. Now, define $\gamma:M_{n,k} \to AM_k$ by $X \mapsto X^TJ_nX$ for all $X \in M_{n,k}$. We can define further, for $f \in \mathcal{P}(AM_k)$, the map
\begin{equation}
	\gamma^*:=f \circ \gamma.
\end{equation}
It follows that $\gamma^* \in \mathrm{Hom}(\mathcal{P}(AM_k),\mathcal{P}(V^{\oplus k})^{\mathrm{Sp}(n)})$ since
\begin{equation}
	\gamma^*(f)(gX)=\gamma^*(f)(X).
\end{equation}
Now,
\begin{equation}
	X \in M_{n,k} \,\Rightarrow\, X^TJ_nX \in AM_k
\end{equation}
and
\begin{equation}
	x_{ij}(X^TJ_nX)=(v_i,J_nv_j) \,\Rightarrow\, \gamma^*(x_{ij})(v_1,\ldots,v_k) = 						\omega(v_i,v_j),
\end{equation}
the contraction of the $i,j$th position with $i<j$ using the $\omega(\cdot\,,\cdot)$.
\begin{theorem}[Polynomial FFT for $\mathrm{Sp}(n)$] 
Let $n \in 2\mathbb{N}$. Then the map $\gamma^*\in \mathrm{Hom} 		
	(\mathcal{P}(AM_k), \mathcal{P}((\mathbb{C}^n)^{\oplus
k})^{\mathrm{Sp}(n)})$ is 					
surjective, and, hence, $\mathcal{P}((\mathbb{C}^n)^{\oplus
k})^{\mathrm{Sp}(n)}$ is 				generated by the
symplectic contractions
	\begin{displaymath}
		\{\omega(v_i,v_j)\,|\,1\leq i< j \leq k \}.
	\end{displaymath}
\end{theorem}

\section{Afterword}

\subsection{Applications}
There are a plethora of applications which use the results of the invariant theory of the classical groups. Here we state just a few. Polynomial and tensor invariants have been used most explicitly in quantum computing \cite{rodriguez}, but have also appeared in works on classical mechanics \cite{morrison}. They also provide a method of describing quantum entanglement \cite{sudbery} and quantum information \cite{brylinski}. However, the results of invariant theory are most often assumed in most all work in special and general relativity where the invariants of $\mathrm{SO}(3,1)$ play a leading role.

\subsection{Comments}
The proof of Theorem \ref{thm:fingen} is due to Hurwitz and follows that presented in \cite{goodman}. The terminology \emph{First Fundamental Theorem} is due to H. Weyl (1946). Also, in his landmark book \emph{The Classical Groups, their Invariants and Representations,} Weyl proved the polynomial form of the FFT using the \emph{Capelli identity} and ``polarization operators''. For a sketch of this method see \cite{fulton}. Also it may be noted that although invariant theory is well-established, ongoing research into novel proofs that may give insight into other fields of mathematics are continuing to be sought after. For example, see \cite{mihailovs1} and \cite{mihailovs2}.

\appendix

\section{The Hilbert Basis Theorem} \label{sec:hilbert}

Since there are two references to this theorem in the body of the paper, it seems necessary to at least state this deep theorem from algebraic geometry. For details and a proof see \cite{goodman}.
\begin{theorem}[Hilbert basis theorem] $\phantom{1}$ \\
	$\phantom{1} \quad$ Let $\mathcal{I} \subseteq \mathcal{P}(V)$ be an ideal. Then 						$\mathcal{I}$ is finitely generated. That is, there is a finite set of polynomials 					$f_1,\ldots,f_d$ in $\mathcal{I}$ so that every $g \in \mathcal{I}$ can be written as
	\begin{equation}
		g=g_1f_1+\cdots+g_df_d
	\end{equation}
	for some choice of $g_1,\ldots,g_d \in \mathcal{P}(V)$.
\end{theorem}

\end{document}